# Comment on "Bayes' Theorem in the 21st Century" by Bradley Efron

Valentin Amrhein[1,2], Tobias Roth[1,2], Fränzi Korner-Nievergelt[3,4*]

[1]Zoological Institute, University of Basel, 4051 Basel, Switzerland; [2]Research Station Petite Camargue Alsacienne, 68300 Saint-Louis, France; [3]oikostat GmbH, 6218 Ettiswil, Switzerland; [4]Swiss Ornithological Institute, 6204 Sempach, Switzerland;
*Correspondence to: fraenzi.korner@vogelwarte.ch

1 August 2013

**In a Perspectives article in *Science* (*1*), Bradley Efron concludes that Bayesian calculations cannot be uncritically accepted when using uninformative priors. We argue that this conclusion is problematic because Efron's example does not use data, hence it is not Bayesian statistics; his priors make little sense and are not uninformative; and using the available data point and an uninformative prior actually leads to a reasonable posterior distribution.**

Efron (*1*) provides four examples of Bayesian analyses, two of which underline the remarkable potential of Bayesian methods. Based on one of the other examples, however, Efron ultimately concludes that Bayesian analyses using uninformative priors cannot be uncritically accepted and should be checked by frequentist methods. While we wholeheartedly agree that statistical results should not be uncritically accepted, we find Efron's example ineffective in showing that Bayesian statistics require more careful checking than any other kind of statistics.

In his example on uninformative priors, Efron uses Bayes' theorem to calculate the probability that twins are identical given that the sonogram shows twin boys. Efron finds this probability to be 2/3 when using an uninformative prior versus 1/2 with an informative prior and thereby concludes that an uninformative prior does not have the desired neutral effects on the output of Bayes' rule. We argue that this example is not only flawed, but completely useless in illustrating Bayesian data analysis because it does not rely on any data. Although there is one data point (a couple is due to be parents of twin boys, and the twins are fraternal), Efron does not use it to update prior knowledge. Instead, Efron combines different pieces of expert knowledge using Bayes' theorem. While certainly an impeccable probability law, Bayes' theorem is a mathematical equation, not a statistical model describing how data may be produced. In essence, Efron uses this equation to show that the value on the left side of the equation changes when a term on the right side is changed, which is trivial and could be shown with any mathematical equation also in a non-Bayesian context. Without new data, our knowledge is by definition determined by prior information; thus, showing that the outcome of a Bayesian analysis with no new data is heavily influenced by the prior would not argue against Bayesian methods. Indeed, without data, Efron's example is not Bayesian statistics and his conclusion about Bayesian statistics based on this example is unjustified.

We also have other more technical issues with Efron's example. Efron interprets the term P(A) on the right side of the equation as the prior on the probability that twins are identical. To make this prior uninformative, it is assigned a value of P(A) = 0.5 (see (*2*); although this is not stated in (*1*)). This uninformative prior is set in contrast to the informative "doctor's prior" of P(A) = 1/3. First, however, the parameter of interest is P(A|B) rather than P(A) according to Efron's study question (see sidebar in (*1*)), thus the focus should be on the appropriate prior for P(A|B). Second, for the uninformative prior, Efron mentions erroneously that he used a uniform distribution between zero and one, which is clearly different from the value of 0.5 that was used. Third, we find it at least debatable whether a prior can be called an uninformative prior if it has a fixed value of 0.5 given without any measurement of uncertainty. For example, if we knew that our chance of winning the next million-dollar jackpot were 50:50, would we really call this uninformative?

If we use the data point together with an uninformative uniform prior on P(A|B) (see sidebar) to determine the probability of identical twins given the twins are two boys, we obtain, with 95% certainty, a probability of between 0.01 and 0.84; if we use a highly informative prior based on information from the doctor and genetics, we obtain a probability of between 0.49 and 0.51. This looks completely reasonable to us, although of course we do not know much more than we knew before because we had only a single data point.

We would very much like to check our calculations using frequentist methods; however, this is impossible because there is only one data point, and frequentist methods generally cannot handle such situations. Although we agree with Efron (*1*) that the choice of the prior is essential, we conclude that his article gives a biased impression of the influence of uninformative priors. In his example using Bayes' theorem, we found no reliable support for his main conclusion that Bayesian calculations cannot be uncritically accepted when using uninformative priors.

---

Study question: What is the probability of identical twins given the twins are two boys?

Data: One pair of twin boys is fraternal.

Data model: x~Binomial($\theta$, n), where $\theta$ is the probability of identical twins given the twins are two boys, x is the number of identical twins in the data, and n is the total number of pairs of twin boys; in our case: x = 0 and n=1.

The posterior distribution p($\theta$|x) is obtained using Bayes' theorem

p($\theta$|x) = p(x|$\theta$)p($\theta$)/p(x)

We use two different priors p($\theta$):

1) Uninformative prior: p($\theta$) = Unif(0,1) = Beta(1,1)

2) Informative prior: using the information from the doctor and from genetics, we are quite sure that $\theta$ must be around 0.5 (*1*). Transforming this information into a statistical distribution yields p($\theta$) = Beta(10000, 10000), which has a mean of 0.5 and a 95% interval of 0.493 – 0.507. [Note that we had to choose the 95% interval arbitrarily because we are not informed about the certainty of the information provided by the doctor and by genetics].

> Given the single parameter Binomial model, x~Binomial($\theta$, n), and the prior p($\theta$) = Beta($\alpha,\beta$), the solution of the Bayesian analysis is given by the posterior distribution p($\theta$|x) = Beta($\alpha$+x,$\beta$+n-x)  [see any Bayesian textbook, e.g. (*3*), p. 34]
>
> The probability of identical twins given the twins are two boys:
>
> 1) Uninformative prior: p($\theta$|x) = Beta(1+x,1+n-x) = Beta(1+0,1+1-0) = Beta(1, 2), which has an expected value of 0.33 and a 95% interval of 0.013 – 0.84.
>
> 2) Informative prior: p($\theta$|x) = Beta(10000+x,10000+n-x) = Beta(10000+0,10000+1-0) = Beta(10000, 10001), which has an expected value of 0.50 and a 95% interval of 0.49 – 0.51.